\documentclass[12pt]{article}

\usepackage{amsmath,amssymb}
\usepackage{graphicx}
\usepackage{tikz}
\usepackage{color,url}

\newtheorem{theorem}{Theorem}[section]
\newtheorem{lemma}[theorem]{Lemma}
\newtheorem{corollary}[theorem]{Corollary}
\newtheorem{proposition}[theorem]{Proposition}

\newtheorem{example}[theorem]{Example}
\newtheorem{conjecture}[theorem]{Conjecture}

\newcommand{\proof}{\noindent{\bf Proof.\ }}
\newcommand{\qed}{\hfill $\square$\medskip}

\tikzstyle{noeud}=[circle,inner sep=2, minimum size =3 pt, line width = 1pt, draw=black, fill=white]
\tikzstyle{noeud2}=[circle,inner sep=1.4, minimum size =3 pt, line width = 1pt, draw=black, fill=white]

\def\cp{\,\square\,}
\DeclareMathOperator {\g} {g}
\DeclareMathOperator {\sg} {sg}

\DeclareMathOperator {\sgc} {sgc}
\DeclareMathOperator {\diam} {diam}

\begin{document}

\title{Strong geodetic cores and Cartesian product graphs}

\author{
    Valentin Gledel $^{a}$
    \and
	Vesna Ir\v si\v c $^{b,c}$
	\and
	Sandi Klav\v zar $^{b,c,d}$
}

\date{\today}

\maketitle
% \vspace{-0.8 cm}
\begin{center}
	$^a$ Univ Lyon, Universit\'e Lyon 1, LIRIS UMR CNRS 5205, F-69621, Lyon, France \\
	\medskip

	$^b$ Institute of Mathematics, Physics and Mechanics, Ljubljana, Slovenia\\
	\medskip

	$^c$ Faculty of Mathematics and Physics, University of Ljubljana, Slovenia\\
% 	{\tt sandi.klavzar@fmf.uni-lj.si}\\
	\medskip
	
	$^d$ Faculty of Natural Sciences and Mathematics, University of Maribor, Slovenia
		
\end{center}

\begin{abstract}
The strong geodetic problem on a graph $G$ is to determine a smallest set of vertices  such that by fixing one shortest path between each pair of its vertices, all vertices of $G$ are covered. To do this as efficiently as possible, strong geodetic cores and related numbers are introduced. Sharp upper and lower bounds on the strong geodetic core number are proved. Using the strong geodetic core number an earlier upper bound on the strong geodetic number of Cartesian products is improved. It is also proved that ${\rm sg}(G \,\square\, K_2) \geq {\rm sg}(G)$ holds for different families of graphs, a result conjectured to be true in general. Counterexamples are constructed demonstrating that the conjecture does not hold in general. 
\end{abstract}


\medskip\noindent
{\bf AMS Subj.\ Class:} 05C12, 05C76

%%%%%%%%%%%%%%%%%%%%%%%%%%%%%%%%%%%%%%%%%%%%%%%%%%%%%%%%%%%%%%%%%%%%%5
%%%%%%%%%%%%%%%%%%%%%%%%%%%%%%%%%%%%%%%%%%%%%%%%%%%%%%%%%%%%%%%%%%%%%5
\section{Introduction}
\label{sec:intro}
%%%%%%%%%%%%%%%%%%%%%%%%%%%%%%%%%%%%%%%%%%%%%%%%%%%%%%%%%%%%%%%%%%%%%5
%%%%%%%%%%%%%%%%%%%%%%%%%%%%%%%%%%%%%%%%%%%%%%%%%%%%%%%%%%%%%%%%%%%%%5

Motivated by social networks applications, the strong geodetic problem was introduced in~\cite{MaKl16a} as follows. Let $G=(V,E)$ be a graph. Given a set $S\subseteq V$, for each pair of vertices $\{x,y\}\subseteq S$, $x\ne y$, let $\widetilde{g}(x,y)$ be a {\em selected fixed} shortest path between $x$ and $y$. Note that this means that $\widetilde{g}(x,y) = \widetilde{g}(y,x)$ for each $x,y \in S$. Set
$$\widetilde{I}(S)=\{\widetilde{g}(x, y) : x, y\in S\}\quad {\rm and}\quad V(\widetilde{I}(S))=\bigcup_{\widetilde{P} \in \widetilde{I}(S)} V(\widetilde{P})\,.$$ 
If $V(\widetilde{I}(S)) = V$ for some $\widetilde{I}(S)$, then the set $S$ is called a {\em strong geodetic set}. For a graph $G$ with just one vertex, we consider the vertex as its unique strong geodetic set. The {\em strong geodetic problem} is to find a smallest strong geodetic set of $G$, such a set will be briefly called an {\em sg-set}. The cardinality of an sg-set is the {\em strong geodetic number} $\sg(G)$ of $G$, that is, if $S$ is an sg-set of $G$, then $\sg(G) = |S|$.  

So far, several aspects of the strong geodetic number have been studied. In the seminal paper~\cite{MaKl16a} the invariant was compared with the isometric path number (see~\cite{clarke-2008, fitzpatrick-1999} for the latter), determined for complete Apollonian networks, and proved that the problem is NP-complete. The strong geodetic number was further investigated on grids and cylinders, and on general Cartesian product graphs in~\cite{products, Klavzar+2017}, respectively. Additional properties, in particular with respect to the diameter, and a solution for balanced complete bipartite graphs were reported in~\cite{bipartite}. The edge version of the problem has also been introduced and investigated in~\cite{MaKl16b}. 

Once $\sg(G)$ is known, we can be faced with sg-sets that have significantly different properties with respect to various applications. Typically we would like to select an sg-set such that covering the graph is as efficient as possible. One way to do it is to involve only those pairs of vertices which cannot be avoided. More formally, if $G=(V(G),E(G))$ is a graph and $S\subseteq V(G)$ is a strong geodetic set, then a set of vertices $X\subseteq S$ is a {\em strong geodetic core for} $S$, if there exists $\widetilde{I}(S)$ such that 
$$\bigcup_{(u,v)\in X\times S} V(\widetilde{g}(u,v)) = V(G)\,.$$
That is, $X$ is a strong geodetic core for $S$ provided that we can cover $V(G)$ without geodesics between the pairs of vertices from $S\setminus X$. The {\em strong geodetic core number} $\sgc(S)$ is the size of a smallest strong geodetic core for $S$, and the {\em strong geodetic core number} of the graph $G$ is  
$$\sgc(G) = \min\{ \sgc(S):\ S \ \text{sg-set}\}\,.$$
If $X$ is a smallest strong geodetic core for $S$, then we will briefly call it an {\em sgc-set}. 

A concept that is in a way dual to the strong geodetic problem is the one of the 
the $k$-path vertex cover~\cite{bresar-2011} where a smallest set of vertices is to be found that hits all the paths of order $k$. This problem has been thoroughly investigated by now, in particular on graph products~\cite{jakovac-2013}. In this paper we in particular investigate the strong geodetic problem on Cartesian products and complete bipartite graphs, and, as a coincidence, the same was done very recently for the $k$-path vertex cover in~\cite{li-2018}. 

The paper is organized as follows.  In the next section we prepare the ground for the rest of the paper: give notation, recall known concepts, introduce generalized geodetic graphs, list a couple of known results and a conjecture, and describe a construction to be used later.  In Section~\ref{sec:examples-and-basics} we prove general upper and lower bounds on the strong geodetic core number and for each of them prove its sharpness. Then, in Section~\ref{sec:complete-bipartite}, we use complete bipartite graphs to obtain several (surprising) properties of the strong geodetic (core) number. In the subsequent section we demonstrate the usefulness of the strong geodetic core number by improving an earlier known   upper bound on the strong geodetic number of Cartesian products. In Section~\ref{sec:valid_conjecture} we show that in many cases the strong geodetic number of a prism over $G$ is at least as large as the strong geodetic number of $G$, a result conjectured in~\cite{products} to be true for all graphs $G$. More generally, similar results are proved for Cartesian products with an arbitrary second factor. In the subsequent section we prove, however, that the conjecture does not hold in general. We conclude the paper with several open problems and directions for further investigation. 

%%%%%%%%%%%%%%%%%%%%%%%%%%%%%%%%%%%%%%%%%%%%%%%%%%%%%%%%%%%%%%%%%%%%%5
%%%%%%%%%%%%%%%%%%%%%%%%%%%%%%%%%%%%%%%%%%%%%%%%%%%%%%%%%%%%%%%%%%%%%5
\section{Preliminaries}
\label{sec:preliminaries}
%%%%%%%%%%%%%%%%%%%%%%%%%%%%%%%%%%%%%%%%%%%%%%%%%%%%%%%%%%%%%%%%%%%%%5
%%%%%%%%%%%%%%%%%%%%%%%%%%%%%%%%%%%%%%%%%%%%%%%%%%%%%%%%%%%%%%%%%%%%%5

The order and the size of a graph $G$ will be denoted with $n(G)$ and $m(G)$, respectively. The distance $d_G(u,v)$ between vertices $u$ and $v$ of a graph $G$ is the usual shortest-path distance. The {\em interval} $I_G(u,v)$ between $u$ and $v$ of a graph $G$ consists of all vertices that lie on some shortest $u,v$-path. If $T$ is a tree, then $\ell(T)$ is the number of leaves of $T$. We will use the notation $[n] = \{1,\ldots, n\}$. 

The {\em Cartesian product} $G\cp H$ of graphs $G$ and $H$ is the graph with vertex set $V(G) \times V(H)$, where the vertices $(g,h)$ and $(g',h')$ are adjacent if either $g=g'$ and $hh'\in E(H)$, or $h=h'$ and $gg'\in E(G)$. If $h\in V(H)$, then a subgraph of $G\cp H$ induced by the set of vertices $\{(x,h) ;\ x\in V(G)\}$ is isomorphic to $G$ and called a {\em $G$-layer}. Analogously $H$-layers are defined. If $X\subseteq V(G\cp H)$, then its  projection $p_G(X)$ on $G$ is  the set $\{g\in V(G):\ (g,h)\in X\ {\rm for\ some}\ h\in V(H)\}$. The projection $p_H(X)$ of $X$ on $H$ is defined analogously. If $X$ is a subgraph of $G\cp H$, then $p_G(X)$ and $p_H(X)$ are the natural projections of $X$ on $G$ and $H$ respectively. Recall that if $P$ is a geodesic in $G\cp H$, then $p_G(P)$ is a geodesic in $G$. For more information on the Cartesian product of graphs see the book~\cite{imrich-2008}. The following conjecture from~\cite{products} was one of the main motivations for the present paper. 

\begin{conjecture}
\label{conj:lower-bound}
If $G$ is a graph with $n(G)\ge 2$, then $\sg(G \cp K_2) \geq \sg(G)$.
\end{conjecture}

In~\cite{products} it was reported, among other related results, that the conjecture holds for all graphs $G$ with $n(G) \leq 7$.

The {\em geodetic number} $\g(G)$ of $G$ is the size of a smallest set $S\subseteq V(G)$ such that $\cup_{\{u,v\}\in \binom{S}{2}} I(u,v) = V(G)$. The  classical {\em geodetic problem} for the graph $G$ is to determine $\g(G)$, see~\cite{bresar-2011b, chartrand-2000, dourado-2010, hansen-2007, pelayo-2013, soloff-2015} for related investigations and in particular~\cite{bresar-2008, jiang-2004} for the geodetic number of Cartesian products. Very recently it was proved in~\cite{bueno-2018} that deciding whether the geodetic number of a graph is at most $k$ is NP-complete for graphs of maximum degree three. This result is appealing since to determine the strong geodetic number is intuitively harder than to determine the geodetic number. 

A graph is called {\em geodetic} if any two vertices are joined by a unique shortest path, cf.~\cite{blokhuis-1988, ore-1962, voblyui-2017}. As observed in~\cite[Lemma 2.1]{MaKl16a}, we have $\sg(G)\ge \g(G)$ for any graph $G$, and if $G$ is geodetic, then $\g(G) = \sg(G)$ holds. Hence we say that $G$ is a {\em generalized geodetic graph} if $\g(G) = \sg(G)$. Consider the split graph $G$ on a clique $K$ of size $m$ and an independent set $I$ of size $n$ and all possible edges between $K$ and $I$. If $\binom{n}{2}\ge m$, then $\g(G) = \sg(G) = n$. This yields a class of generalized geodetic graphs that are not geodetic. 

The (open) neighborhood of a vertex $v$ is denoted with $N(v)$. A vertex $v$ of a graph $G$ is {\em simplicial} if $N(v)$ induces a complete graph. 
For later use we state the following known fact. 

\begin{lemma}
\label{lem:simplicial}
If $S$ is a strong geodetic set of $G$ and $v$ is a simplicial vertex of $G$, then $v\in S$.
\end{lemma}

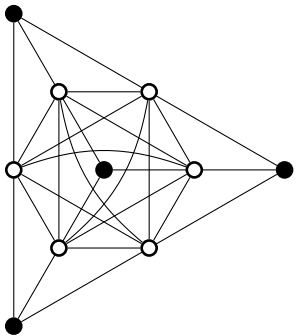
\begin{figure}[!ht]
	\begin{center}
		\begin{tikzpicture}[scale=0.6]
		
		\node[noeud, fill=black] (1) at (1*360/3:4 cm) {};
		\node[noeud, fill=black] (2) at (2*360/3:4 cm) {};
		\node[noeud, fill=black] (3) at (3*360/3:4 cm) {};
		
		\path (1) edge (2);
		\path (2) edge (3);
		\path (1) edge (3);
		
		\node[noeud, fill=black] (0) at (0, 0) {};
		
		\draw (0) edge (1); 
		\path (0) edge (2); 
		\path (0) edge (3); 
		
		\node[noeud] (a) at (1.5*360/3:2 cm) {};
		\node[noeud] (b) at (2.5*360/3:2 cm) {};
		\node[noeud] (c) at (0.5*360/3:2 cm) {};
		\node[noeud] (a') at (1*360/3:2 cm) {};
		\node[noeud] (b') at (2*360/3:2 cm) {};
		\node[noeud] (c') at (3*360/3:2 cm) {};
		
		\draw (a) edge (b);
		\draw (b) edge (c);
		\draw (c) edge (a);
		\draw (a') edge (b');
		\draw (b') edge (c');
		\draw (c') edge (a');
		\draw (a) edge (a');
		\draw (a') edge(c);
		\draw (a) edge (b');
		\draw (b) edge (b');
		\draw (b) edge (c');
		\draw (c) edge (c');

		\draw [bend left=20,-] (b) to (a');
		\draw [bend left=20,-] (c) to (b');
		\draw [bend left=20,-] (a) to (c');
		
		\end{tikzpicture}
		\caption{The graph $\widehat{S(K_4)}$ with its strong geodetic set.}
		\label{fig:sgc_max}
	\end{center}
\end{figure}

As an application of Lemma~\ref{lem:simplicial} we consider the following class of graphs that will be applied later. Let $S(G)$ be the {\em subdivision graph} of a graph $G$, that is, the graph obtained from $G$ by subdividing each of its edges exactly once. In addition, let $\widehat{S(G)}$ be the graph obtained from $S(G)$ by adding an edge between different vertices $u, v\in V(S(G))\setminus V(G)$. See Fig.~\ref{fig:sgc_max} where the graph $\widehat{S(K_4)}$ is drawn.

\begin{proposition}
\label{prop:S(G)-and-more}
If $G$ is a connected graph, then 
$$\sg(\widehat{S(G)}) = n(G)\quad {\rm and}\quad \sgc(\widehat{S(G)}) \ge \min\{k:\  \sum_{i=1}^k (n-i) \ge m(G)\}\,.$$ 
\end{proposition}

\proof
The vertices of $G$ are simplicial in $\widehat{S(G)}$, hence in view of Lemma~\ref{lem:simplicial} we have $\sg(\widehat{S(G)}) \ge n(G)$. On the other hand, these vertices also form a strong geodetic set, hence the first assertion. It also follows that $V(G)$ is a unique sg-set of $\widehat{S(G)}$. 

Note that geodesics between vertices of $V(G)$ are either of length $2$ or $3$. Hence such a geodesic from a vertex $v \in V(G)$ to another vertex of $V(G)$ covers one or two vertices from $V(\widehat{S(G)}) - V(G)$. Moreover, if it covers two, then one is a vertex obtained by subdividing an edge $vx$, and this vertex is also covered with the (unique) $v,x$-geodesic. It follows that geodesics from $v \in V(G)$ to the vertices of $V(G)\setminus \{v\}$ cover at most $n-1$  vertices in $V(\widehat{S(G)}) - V(G)$. 
After $v$ is added to a core, by the same reasoning we infer that the geodesics between $u \in V(G) - \{v\}$ and the remaining vertices cover at most $n-2$ different vertices. Repeating the same reasoning until all vertices in $V(\widehat{S(G)}) - V(G)$ are covered we obtain $\sgc(\widehat{S(G)}) \ge \min\{k:\  \sum_{i=1}^k (n-i) \ge m(G)\}$. 
\qed

In Section~\ref{sec:complete-bipartite} we will demonstrate that the inequality from Proposition~\ref{prop:S(G)-and-more} can be strict.

%%%%%%%%%%%%%%%%%%%%%%%%%%%%%%%%%%%%%%%%%%%%%%%%%%%%%%%%%%%%%%%%%%%%%5
%%%%%%%%%%%%%%%%%%%%%%%%%%%%%%%%%%%%%%%%%%%%%%%%%%%%%%%%%%%%%%%%%%%%%5
\section{Bounds on the strong geodetic core number}
\label{sec:examples-and-basics}
%%%%%%%%%%%%%%%%%%%%%%%%%%%%%%%%%%%%%%%%%%%%%%%%%%%%%%%%%%%%%%%%%%%%%5
%%%%%%%%%%%%%%%%%%%%%%%%%%%%%%%%%%%%%%%%%%%%%%%%%%%%%%%%%%%%%%%%%%%%%5

If $n\ge 1$, then $\sgc(K_n) = 1$. The same holds for cycles, graphs on at most five vertices, and trees. We state the latter fact as a lemma for later use. 

\begin{lemma}
\label{lem:trees:esc=1}
If $T$ is a tree, then $\sgc(T) = 1$. 
\end{lemma}

\proof
The set $L$ of leaves of $T$ forms a unique sg-set of $T$. If $u\in L$, then consider the BFS-tree rooted at $u$ to see that the geodesics between $u$ and all the other vertices of $L$ cover $V(T)$. Hence $\sgc(T) = 1$. 
\qed

As said, $\sgc(K_n) = 1$. For all the other graphs we have the following bounds. 

\begin{theorem}
\label{thm:general-bounds} 
If $G$ is a graph with $n=n(G)$, $s=\sg(G)$, and $d = \diam(G)\ge 2$, then 
$$1\le \left\lceil s - \frac{1+\sqrt{(2s-1)^2-\frac{8(n-s)}{d-1}}}{2} \right\rceil  \le \sgc(G) \le \min\{s - 1, n - s\}\,.$$
\end{theorem}

\proof
The bound $\sgc(G)\le \sg(G)-1$ holds true because if $S$ is a strong geodetic set and $u\in S$, then setting $X = S\setminus \{u\}$ we have $\cup_{(u,v)\in X\times S}\ \widetilde{g}(u,v) = V(G)$. Since $G$ is not complete, in a (smallest) strong geodetic core $X\subseteq S$ we need at most one vertex for each vertex from $V(G)\setminus S$, hence $\sgc(G)\le n(G)-\sg(G)$.  This proves the upper bounds. 

For the lower bound we claim that the following inequality holds:
\begin{equation}
	\label{eqn:sgc_lower}
	\left(\sgc(G)(\sg(G)-\sgc(G))+\binom{\sgc(G)}{2}\right)(\diam(G)-1) \ge n(G) - \sg(G)\,.
\end{equation}
Let $S$ and $C$ be respectively an sg-set of $G$ and its core such that $|C| = \sgc(G)$. The product $\sgc(G)(\sg(G)-\sgc(G))$ is the number of paths between the vertices from $C$ and the vertices from $S\setminus C$, while $\binom{\sgc(G)}{2}$ is the number of paths between two vertices $C$. Each of these paths cover at most $\diam(G)-1$ vertices of $V(G) \setminus S$. Since only the paths between $C$ and $S$ are needed to cover $V(G)$, these paths must cover every vertex of $V(G)\setminus S$. This proves the inequality. After solving it, the second lower bound is obtained. 

Recall from~\cite{MaKl16a} that $\sg(G) < n(G)$ holds for all non complete graphs. It follows that $8(n(G)-\sg(G))/(\diam(G)-1) > 0$. Hence the expression from the theorem under the ceiling is strictly positive which proves the first lower bound. 
\qed

We next show that all the bounds of Theorem~\ref{thm:general-bounds} are sharp. The sharpness of $\sgc(G) \ge 1$ follows from Lemma~\ref{lem:trees:esc=1}. For the other bounds we have: 

\begin{theorem}
\label{thm:bounds-are-sharp} 
\begin{enumerate}
\item[(i)] For every $k\ge 2$ there exists a graph $G$ with $\sgc(G) = k$ and $\sg(G) = k+1$. 
\item[(ii)] For every $n\ge 3$ there exists a graph $G$ of order $n$ with $\sgc(G) = \left\lfloor \frac{n}{3} \right\rfloor$ and $\sg(G) = n - \left\lfloor \frac{n}{3} \right\rfloor$.  
\item[(iii)] For every $k\ge 1$, $d\ge 2$, $s\ge 0$, there exists a graph $G$ such that 
$$\sgc(G) = k =  \left\lceil \sg(G) - \frac{1+\sqrt{(2\sg(G)-1)^2-\frac{8[n(G)-\sg(G)]}{{\rm diam}(G)-1}}}{2} \right\rceil\,,$$
where $\sg(G) = k + s$ and ${\rm diam}(G) = d$. 
\end{enumerate}
\end{theorem}

\proof
(i) Let $G_k = \widehat{S(K_k)}$, $k\ge 3$. (Recall that $G_4$ is drawn in Fig.~\ref{fig:sgc_max}.) Then from Proposition~\ref{prop:S(G)-and-more} we deduce that $\sgc(G_{k+1}) \ge k$ and $\sg(G_{k+1}) = k+1$ for every $k \geq 2$.  Since for any graph $G$ we have $\sgc(G)\le \sg(G) - 1$, we also have $\sgc(G_{k+1}) \le k$ and we are done. 

(ii) If $n$ is even, then consider the cocktail-party graph $G = K_{2,\ldots,2}$ of order $n$. If $n$ is odd, then let $G = K_{2,\ldots,2, 1}$ be of order $n$. We claim that $\sg(G) = n - \left\lfloor n/3\right\rfloor$ and $\sgc(G) = \left\lfloor n/3\right\rfloor$. Let $S$ be a minimum strong geodetic set of $G$. Then vertices $u,v\in S$ can cover an additional vertex if and only if $u$ and $v$ are not adjacent, that is, are in the same part of the cocktail party graph. To minimize $|S|$, the set $S$ contains as many non-adjacent pairs of vertices as possible, that is, $\left\lfloor n/3\right\rfloor$ such pairs. These pairs cover $3\left\lfloor n/3\right\rfloor$ vertices, while the remaining ($0$, $1$, or $2$) vertices must be covered by themselves. Therefore, 
$$\sg(G) = 2\left\lfloor \frac{n}{3}\right\rfloor + (n - 3\left\lfloor \frac{n}{3}\right\rfloor) = n - \left\lfloor \frac{n}{3}\right\rfloor\,.$$
To see that $\sgc(G) = \left\lfloor n/3\right\rfloor$, consider a minimum strong geodetic set $S$ of $G$. Then, by the above, each vertex $x\in V(G)\setminus S$ is covered by a unique geodesic between non-adjacent vertices $u,v\in S$. It follows that a strong geodetic core must contain one of $u$ and $v$. Selecting one vertex from each pair of non-adjacent vertices from $S$ yields a strong geodetic core, hence $\sgc(G) = \left\lfloor n/3\right\rfloor$ indeed holds. 

\begin{figure}[!ht]
	\begin{center}
		\begin{tikzpicture}[scale=0.6]
		
		\node[noeud, fill=black] (1) at (1*360/3:4 cm) {};
		\node[noeud, fill=black] (2) at (2*360/3:4 cm) {};
		\node[noeud, fill=black] (3) at (3*360/3:4 cm) {};
		
		\path (1) edge (2);
		\path (2) edge (3);
		\path (1) edge (3);
		
		\node[noeud, fill=black] (0) at (2.5*360/3:1 cm) {};
		\node[noeud, fill=black] (0') at (3.5*360/3:1 cm) {};
		
		\path (0) edge (0');
		
		\draw (0) edge (1); 
		\path (0) edge (2); 
		\path (0) edge (3);
		\draw (0') edge (1); 
		\path (0') edge (2); 
		\path (0') edge (3); 
		
		\node[noeud, fill=gray] (a) at (1.5*360/3:2 cm) {};
		\node[noeud, fill=gray] (b) at (2.5*360/3:2 cm) {};
		\node[noeud, fill=gray] (c) at (0.5*360/3:2 cm) {};
		\node[noeud, fill=gray] (a') at (1*360/3:2 cm) {};
		\node[noeud, fill=gray] (a'') at (0.88*360/3:2 cm) {};
		\node[noeud, fill=gray] (b') at (2*360/3:2 cm) {};
		\node[noeud, fill=gray] (b'') at (2.12*360/3:2 cm) {};
		\node[noeud, fill=gray] (c') at (3.12*360/3:2 cm) {};
		\node[noeud, fill=gray] (c'') at (2.88*360/3:2 cm) {};
		
		\node[noeud2] (x1) at (2.75*360/3:2.3 cm) {};
		\node[noeud2] (x2) at (2.25*360/3:2.3 cm) {};
		\node[noeud2] (x3) at (3.25*360/3:2.3 cm) {};
		\node[noeud2] (x4) at (3.75*360/3:2.3 cm) {};
		\node[noeud2] (x5) at (1.25*360/3:2.3 cm) {};
		\node[noeud2] (x6) at (1.75*360/3:2.3 cm) {};
		
		\node[noeud2] (y1) at (2*360/3:2.6 cm) {};
		\node[noeud2] (y2) at (2.07*360/3:2.6 cm) {};
		\node[noeud2] (y3) at (1*360/3:2.6 cm) {};
		\node[noeud2] (y4) at (0.93*360/3:2.6 cm) {};
		\node[noeud2] (y5) at (3.065*360/3:2.6 cm) {};
		\node[noeud2] (y6) at (2.935*360/3:2.6 cm) {};
		
		\node[noeud2] (z1) at (1*360/3:1.4 cm) {};
		\node[noeud2] (z2) at (0.75*360/3:1.4 cm) {};
		\node[noeud2] (z3) at (2*360/3:1.4 cm) {};
		\node[noeud2] (z4) at (2.25*360/3:1.4 cm) {};
		\node[noeud2] (z5) at (3.25*360/3:1.4 cm) {};
		\node[noeud2] (z6) at (2.75*360/3:1.4 cm) {};
		
		\end{tikzpicture}
		\caption{The representation of the graph $H_{3,2,4}$ with its strong geodetic set in black and the set $\mathcal{V}$ in gray.}
		\label{fig:sgc_min}
	\end{center}
\end{figure}
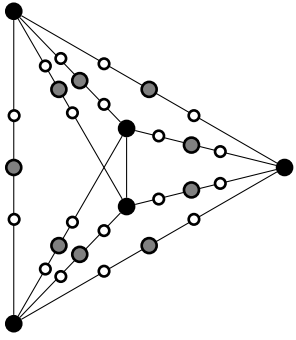

(iii) Consider a graph obtained from the join of the graphs $K_k$ and $K_s$ by subdividing every edge $d-1$ times, except the edges in $K_s$, which remain unsubdivided. Set $V_{k,s} = V(K_k) \cup V(K_s)$, and let the vertices of the graph be 
$$V_{k,s} \cup \{ e^i \; ; \; e \in E(K_k), i \in [d-1] \} \cup \{ f_{u,v}^i \; ; \; u \in V(K_k), v \in V(K_s), i \in [d-1] \}.$$
Define a set 
\begin{align*}
	\mathcal{V} = & \{ e^i \; ; \; e \in E(K_k), i \in \{ \lfloor d/2 \rfloor, \lceil d/2 \rceil \} \} \cup \\
	        \cup  & \{ f_{u,v}^i \; ; \; u \in V(K_k), v \in V(K_s), i \in \{ \lfloor d/2 \rfloor, \lceil d/2 \rceil \} \},
\end{align*}
and add all missing edges between vertices in $\mathcal{V}$ (hence $\mathcal{V}$ forms a clique) to obtain the graph $H_{k,s,d}$ (cf.~Fig.~\ref{fig:sgc_min} for $H_{3,2,4}$). Clearly, $n(H_{k,s,d}) = k + s + (ks + \binom{k}{2}) (d-1)$. Moreover, it is straighforward to see that $\diam(H_{k,s,d}) = d$. We claim that $\sg(H_{k,s,d}) = k + s$, and $\sgc(H_{k,s,d}) = k$. Note that this implies that the graph $H_{k,s,d}$ attains the equality in~\eqref{eqn:sgc_lower}, and thus also in $$k =  \left\lceil \sg(G) - \frac{1+\sqrt{(2\sg(G)-1)^2-\frac{8[n(G)-\sg(G)]}{{\rm diam}(G)-1}}}{2} \right\rceil\,.$$

To show this, consider the following subsets defined for all $v \in V_{k,s}$:
$$M_v = \{v\} \cup \{ x:\ x \in V(H_{k,s,d}) - \mathcal{V}, d(v,x) < d(v',x)\ {\rm for\ every}\ v' \in V_{k,s}\setminus \{v\}\}\,.$$
It follows from the geodetic structure of the graph that if none of the vertices in some $M_v$ lies in a strong geodetic set, then no geodesic covers $v$. Moreover, if $v$ is not in a strong geodetic set, then for each edge incident with $v$, at least one vertex from $M_v$ must belong to the strong geodetic set. Note that if there is just one such edge, then $v$ is simplicial. Thus each minimal strong geodetic set contains vertices $V_{k,s}$. As this set is indeed a strong geodetic set, it follows that $\sg(H_{k,s,d}) = k+s$, where $V_{k,s}$ is a unique sg-set.

As for the graph $G_k$ defined in part (i) of the proof, it holds that at least $k-1$ vertices are in the core of $H_{k,s,d}$. But as this is clearly not enough and $V(K_k)$ is a strong geodetic core, we have $\sgc(H_{k,s,d}) = k$.  
\qed

If a graph $G$ attains the equality $\sgc(G) = \sg(G) - 1$, then each of the $\binom{\sg(G)}{2}$ geodesics must cover at least one private vertex. Thus such a graph must have at least $\sg(G) + \binom{\sg(G)}{2}$ vertices. As for the graph $G_{k+1}$ from the proof of Theorem~\ref{thm:bounds-are-sharp}(i) we have $|V(G_{k+1})| = k+1 + \binom{k+1}{2}$, these graphs are smallest possible examples for this situation.

%%%%%%%%%%%%%%%%%%%%%%%%%%%%%%%%%%%%%%%%%%%%%%%%%%%%%%%%%%%%%%%%%%%%%5
%%%%%%%%%%%%%%%%%%%%%%%%%%%%%%%%%%%%%%%%%%%%%%%%%%%%%%%%%%%%%%%%%%%%%5
\section{On the strong geodetic problem in complete bipartite graphs}
\label{sec:complete-bipartite}
%%%%%%%%%%%%%%%%%%%%%%%%%%%%%%%%%%%%%%%%%%%%%%%%%%%%%%%%%%%%%%%%%%%%%5
%%%%%%%%%%%%%%%%%%%%%%%%%%%%%%%%%%%%%%%%%%%%%%%%%%%%%%%%%%%%%%%%%%%%%5

The simple class of complete bipartite graphs turned out as a big challenge for the strong geodetic problem. As already mentioned, $\sg(K_{n,n})$ was determined in~\cite{bipartite}, while such a result for all complete bipartite graphs is yet to be found. In this section we show further appealing properties of complete bipartite graphs with respect to the strong geodetic problem. We first observe that if $n\ge 1$, then 
\begin{equation}
\label{eq:complete-bipartite-special}
\sg\big(K_{n,\binom{n}{2}}\big) = n\,.
\end{equation}
Indeed, the bipartition set of order $n$ is a strong geodetic set as the pairs of its vertices can take care for exactly $\binom{n}{2}$ vertices of the other bipartition set. On the other hand, if $\sg\big(K_{n,\binom{n}{2}}\big) \leq n-1$, then geodesic between these vertices would cover at most $\binom{n-1}{2}$ other vertices. Thus at most $(n-1) + \binom{n-1}{2} = \binom{n}{2} < |V\big(K_{n,\binom{n}{2}}\big)|$ vertices would be covered. This proves~\eqref{eq:complete-bipartite-special}. 

We next show that with the help of the complete bipartite graphs $K_{n,\binom{n}{2}}$, $n\ge 4$ we can show that the inequality of Proposition~\ref{prop:S(G)-and-more} cannot be turned into equality. 

\begin{proposition}
\label{prop:non-equality}
If $n\ge 4$ and $G_n = K_{n,\binom{n}{2}}$, then $\sgc(\widehat{S(G_n})) = n$.
\end{proposition}

\proof
Denote by $(X,Y)$ the bipartition of $G_n$, where $|X| = n$, $|Y| = \binom{n}{2}$.  

First observe that $\min\{k:\  \sum_{i=1}^k (n-i) \ge m(G_n)\} = n-1$, hence Proposition~\ref{prop:S(G)-and-more} yields $\sgc(\widehat{S(G_n)}) \geq n-1$. Suppose the equality is attained. As observed in the proof of Proposition~\ref{prop:S(G)-and-more}, $S = V(G_n)$ is a unique sg-set. Suppose that $C$ is its core of size $n-1$ and let $k = |C \cap Y|$. 

Vertices with the neighborhood in $S - C$ can be covered only by geodesics between $X \cap C$ and $X - C$, or between $Y \cap C$ and $Y - C$. Hence,
$$(n - (n-1-k)) \left( \binom{n}{2} - k\right)  \leq (n-1-k) (n - (n-1-k)) + k \left( \binom{n}{2} - k\right)\, ,$$
which simplifies to 
$$\binom{n}{2} - k - (n-k-1)(k+1) \leq 0 \, .$$
But the quadratic equation in $k$ on the left side has the minimum value $\frac{(n-1) (n-3)}{4}$, which is strictly positive for $n \geq 4$. Hence, $\sgc(\widehat{S(G_n)}) > n-1$. Note that on the other hand, $\sgc(\widehat{S(G_n)}) \leq n$, as $X$ is a strong geodetic core. We conclude that $\sgc(\widehat{S(G_n)}) = n > n-1$. 
\qed

As we have seen in the above proof that  $\min\{k:\  \sum_{i=1}^k (n-i) \ge m(G_n)\} = n-1$, Proposition~\ref{prop:non-equality} thus indeed shows that the inequality in Proposition~\ref{prop:S(G)-and-more} is strict for $K_{n,\binom{n}{2}}$. 

\medskip
To conclude the section, we use complete bipartite graphs again to present an interesting phenomenon on the strong geodetic core number. If we have strong geodetic sets $S \subseteq T$ of a graph $G$, then clearly $\sgc(S) \geq \sgc(T)$. On the other hand, as the following ``7/11-example" shows,  there exist strong geodetic sets $S$ and $T$ with $|S| \leq |T|$ such that $\sgc(S) < \sgc(T)$. (By the above, in such a case $T$ cannot be a superset of $S$.) 

\begin{example} {\rm (7/11 example)}
	Consider the complete bipartite graph $K_{7,11}$ with vertex set $\{x_1, \ldots, x_{11}\} \cup \{y_1, \ldots, y_7\}$. Notice that the sets $S = \{x_1, \ldots, x_7\}$ and $T = \{ x_1, \ldots, x_5 \} \cup \{ y_1, \ldots, y_3 \}$ are strong geodetic sets. In fact, by a simple case analysis we can prove that $S$ is the unique optimal strong geodetic set. Notice that $\sgc(S) = 2$ (an appropriate set $X$ is $\{ x_1, x_2 \}$), and $\sgc(T) = 4$ (an appropriate set is $\{ x_1, x_2, x_3, y_1 \}$). 
\end{example}

%%%%%%%%%%%%%%%%%%%%%%%%%%%%%%%%%%%%%%%%%%%%%%%%%%%%%%%%%%%%%%%%%%%%%5
%%%%%%%%%%%%%%%%%%%%%%%%%%%%%%%%%%%%%%%%%%%%%%%%%%%%%%%%%%%%%%%%%%%%%5
\section{An application to Cartesian products}
\label{sec:products-upper-bounds}
%%%%%%%%%%%%%%%%%%%%%%%%%%%%%%%%%%%%%%%%%%%%%%%%%%%%%%%%%%%%%%%%%%%%%5
%%%%%%%%%%%%%%%%%%%%%%%%%%%%%%%%%%%%%%%%%%%%%%%%%%%%%%%%%%%%%%%%%%%%%5

In~\cite{products} the strong geodetic number was studied on Cartesian products. In particular, it was proved (see~\cite[Theorem 2.1]{products}) that if $G$ and $H$ are graphs, then
\begin{equation}
\label{eq:old-upper-Cartesian}
\sg(G \cp H) \leq \min\{ \sg(H) n(G) - \sg(G) + 1, \sg(G) n(H) - \sg(H) + 1 \}\,.
\end{equation}
We can use the concept of the strong geodetic core to improve this result as follows. 

\begin{theorem}
	\label{thm:upper_bound_sgc}
	If $G$ and $H$ are graphs, then
	$$\sg(G \cp H) \leq \min\{ \sgc(H)(n(G)-1) + \sg(H), \sgc(G)(n(H)-1) +\sg(G) \}\,.$$
\end{theorem}

\proof
By commutativity of the Cartesian product operation, it suffices to prove that  $\sg(G \cp H) \leq \sgc(H)(n(G)-1) + \sg(H)$. 

Let $S_H$ be an sg-set of $H$ such that $\sgc(S_H)=\sgc(H)$, $C_H$ be a strong geodesic core of $S_H$ and $\widetilde{I}(S_H)$ fixed geodesics in $H$. Set $l = |S_H| = \sg(H)$, $m = |C_H| = \sgc(H)$ and $S_H = \{ h_0, h_1, \ldots, h_{l-1}\}$, with $C_H = \{ h_0, h_1, \ldots, h_{m-1}\}$ .
Denote with $Q_{i,j}$ the $h_i, h_j$-geodesic from $\widetilde{I}(S_H)$ for all $i,j \in [0,l-1]$. Fix a vertex $g_0 \in V(G)$, and shortest paths $P_g$ from $g$ to $g_0$ in $G$ for all $g \in V(G) - \{ g_0 \}$. 

Define $T = ((V(G)\setminus\{g_0\}) \times C_H) \cup (\{g_0\} \times S_H))$. Clearly, $|T| = \sgc(H)(n(G)-1) + \sg(H)$. We claim that $T$ is a strong geodetic set of $G \cp H$. To show it, we first fix geodesics in $H$-layers between vertices from $T$ in the same way as they are fixed in $\widetilde{I}(S_H)$. The only (possibly) uncovered vertices are the ones lying in $H$-layers $^{g}H$ for $g \in V(G) - \{ g_0 \}$ that lie on paths $^{g}Q_{i,j}$ for $i \in [0, l-1]$ and $j \in [m,l-1]$. Since the vertices $\{ h_m, \ldots, h_{l-1}\}$ are not in $C_H$, the paths $Q_{i,j}$ with $i,j \in [m, l-1]$ are not needed to cover $H$ and the corresponding paths $^{g}Q_{i,j}$ are not needed to cover the layer $^{g}H$. Only the vertices on the paths $^{g}Q_{i,j}$ for $i \in [0, m-1]$ and $j \in [m,l-1]$ still need to be covered.

To cover them we fix $(g,h_i), (g_0, h_j)$-geodesics as paths $^{g}Q_{i,j}$ joined with $P_g^{h_j}$ for all $g \in V(G) \setminus \{ g_0 \}, i \in [0, m-1]$ and $j \in [m, l-1]$. In this way all the vertices of $G \cp H$ are covered, hence $\sg(G \cp H) \leq |T|$.
\qed

Quite often $\sgc(G)$ is significantly smaller that $\sg(G)$; as already noticed, this happens for complete graphs and for trees. In such cases Theorem~\ref{thm:upper_bound_sgc} gives much better bounds than~\eqref{eq:old-upper-Cartesian}. Consider for instance $K_n\cp K_n$, $n\ge 1$. Then~\eqref{eq:old-upper-Cartesian} yields $\sg(K_n \cp K_n) \leq n^2 - n + 1$, while Theorem~\ref{thm:upper_bound_sgc} gives $\sg(K_n \cp K_n) \leq 2n -1$, which is the exact value, see~\cite[Theorem 3.3]{products}. But even if $\sgc(G)$ is large, Theorem~\ref{thm:upper_bound_sgc} is generally better than~\eqref{eq:old-upper-Cartesian}. Indeed, since $\sgc(G) \le \sg(G) - 1$, Theorem~\ref{thm:upper_bound_sgc} implies the following result stronger than~\eqref{eq:old-upper-Cartesian}:

\begin{corollary}
	\label{prop:upper_bound}
	If $G$ and $H$ are graphs, then
	$$\sg(G \cp H) \leq \min\{ \sg(H)n(G) - n(G) + 1, \sg(G)n(H) - n(H) + 1 \}\,.$$
\end{corollary}

We conclude the section with another family for which the bound in Theorem~\ref{thm:upper_bound_sgc} is attained. 

\begin{proposition}
\label{prop:P_n-P_3}
If $n\ge 3$, then $\sg(P_n\cp P_3) = 4$. 
\end{proposition}

\proof
Theorem~\ref{thm:upper_bound_sgc} yields $\sg(P_n\cp P_3) \le 4$. To prove the other inequality, assume on the contrary that $\sg(P_n\cp P_3) = 3$. Let $S = \{x,y,z\}$ be an sg-set of $P_n\cp P_3$. Note that $S$ must contain at least one vertex in both  $P_3$-layers above the leaves of $P_n$. Let $x, y$ be these two vertices, then $d(x,y)\in \{n-1, n,n+1\}$. If $d(x,y) = n+1$, then $d(x,z) + d(z,y) = n+1$, if $d(x,y) = n$, then $d(x,z) + d(z,y)\le n+2$, and if $d(x,y) = n-1$, then $d(x,z) + d(z,y)\le n+3$. In any case, the union of an $x,y$-geodesic, an $x,z$-geodesic, and a $y,z$-geodesic cover at most $2n+5$ vertices, where $x,y,z$ are counted twice, hence we can cover at most $2n+2 < 3n$ vertices. 
\qed

%%%%%%%%%%%%%%%%%%%%%%%%%%%%%%%%%%%%%%%%%%%%%%%%%%%%%%%%%%%%%%%%%%%%%5
%%%%%%%%%%%%%%%%%%%%%%%%%%%%%%%%%%%%%%%%%%%%%%%%%%%%%%%%%%%%%%%%%%%%%5
\section{Valid cases for Conjecture~\ref{conj:lower-bound}}
\label{sec:valid_conjecture}
%%%%%%%%%%%%%%%%%%%%%%%%%%%%%%%%%%%%%%%%%%%%%%%%%%%%%%%%%%%%%%%%%%%%%5
%%%%%%%%%%%%%%%%%%%%%%%%%%%%%%%%%%%%%%%%%%%%%%%%%%%%%%%%%%%%%%%%%%%%%5

In this section we show that Conjecture~\ref{conj:lower-bound} holds for several large families of graphs. Actually the more general problem~\cite[Problem 3.5]{products} asks whether $\sg(G \cp H) \geq \max\{ \sg(G), \sg(H) \}$ holds for arbitrary graphs $G$ and $H$. Here we prove several cases in which the answer to this question is positive.  

\begin{lemma}
\label{lem:projection}
If $G$ and $H$ are graphs and $S$ is a strong geodetic set of $G \cp H$, then $p_G(\widetilde{I}(S))$ covers every vertex of $G$ and $p_G(S)$ is a geodetic set of $G$.
\end{lemma}

\proof
Let $S_G = p_G(S)$ and $\widetilde{I}(S_G) = p_G(\widetilde{I}(S))$. If $(u,v)\in V(G \cp H)$, then there exists a path $Q$ in $\widetilde{I}(S)$ that covers $(u,v)$. Then $p_G(Q)$ covers $u$ and, since it is in $\widetilde{I}(S_G)$, the set of paths $\widetilde{I}(S_G)$ covers all vertices of $G$. Since these projected paths connect vertices from $S_G$, we conclude that $S_G$ is a geodetic set of $G$.
\qed

\begin{proposition}
\label{prp:gen_geodesic}
If $G$ is a generalized geodetic graph, then $\sg(G\cp H)\ge \sg (G)$ holds for every graph $H$. 
\end{proposition}

\proof
Let $S$ be an sg-set of $G\cp H$ and let $S_G = p_G(S)$. By Lemma~\ref{lem:projection}, $S_G$ is a geodetic set of $G$, hence $|S_G| \ge \g(G)$. Since $G$ is a generalized geodetic graph, $\g(G)=\sg(G)$. As $S_G$ is the projection of $S$ we have also have $|S_G| \le |S|$. Putting these estimates together we get
$$\sg(G) = \g(G) \le |S_G| \le |S| = \sg(G \cp H)\,,$$
hence the result. 
\qed

\begin{corollary}
\label{cor:gen_geodetic}
If $G$ and $H$ are generalized geodetic graphs, then $\sg(G \cp H) \geq \max\{ \sg(G), \sg(H) \}$. 
\end{corollary}

Note that Corollary~\ref{cor:gen_geodetic} in particular confirms Conjecture~\ref{conj:lower-bound} for prisms over generalized geodetic graphs.  

When equality holds in Proposition~\ref{prp:gen_geodesic} we can say more about the structure of the corresponding sg-set. 

\begin{proposition}
\label{prp:eq_gengeo}
Let $G$ be a generalized geodetic graph and $H$ a graph such that $\sg(G \cp H) = \sg(G)$. If $S$ is an sg-set of $G \cp H$, then $|p_G(S)| = |S|$ and $p_G(S)$ is an sg-set of $G$. 
\end{proposition}

\proof
By Lemma~\ref{lem:projection} we know that $p_G(S)$ is a geodetic set of $G$ so $|p_G(S)| \ge \g(G) = \sg(G) = |S|$. Since $p_G(S)$ is the projection of $S$ on $G$ we have the equality $|p_G(S)| = |S|$. Let $\widetilde{I}(S)$ be a set of paths corresponding to the sg-set $S$ and let $\widetilde{I}(S_G) = p_G(\widetilde{I}(S))$. Since $|p_G(S)| = |S|$, the projection is a bijection and between two vertices of $p_G(S)$ there is only one path in $\widetilde{I}(S_G)$. Using Lemma~\ref{lem:projection} again we infer that $\widetilde{I}(S_G)$ covers all vertices of $G$ so $p_G(S)$ is a strong geodetic set of $G$. Since $|p_G(S)| = \sg(G)$, it is an sg-set of $G$.
\qed

Proposition~\ref{prp:eq_gengeo} can in turn be used to find cases in which the bound cannot be reached. To state such a result, we need the following concepts. A subgraph $H$ of a graph $G$ is {\em convex} if for every $u,v\in V(H)$, every shortest $u,v$-path in $G$ lies completely in $H$. If a graph $G$ admits a partition of its vertex set into $p$ non-empty sets $V_1,\ldots, V_p$ such that every $V_i$, $i\in [p]$, induces a convex subgraph, then $G$ is said to admit a {\em convex $p$-partition}, see~\cite{artigas-2011}. The problem whether a graph admits such a partition is NP-complete in general, but polynomial (for every fixed $p\ge 1$) for bipartite graphs~\cite{grippo-2016}.

\begin{theorem}
\label{thm:sgc_sup_sg2}
If $G$ is a generalized geodetic graph with $\sgc(G) >\frac{\sg(G)}{2}$ and $H$ admits a convex $2$-partition, then $\sg(G \cp H) > \sg (G)$.
\end{theorem}

\proof
Assume on the contrary that $\sg(G \cp H) \le \sg (G)$. Let $S$ be an sg-set of $G \cp H$ and $S_G = p_G(S)$. By Proposition~\ref{prp:gen_geodesic} we have $\sg(G \cp H) = \sg (G)$  and hence Proposition~\ref{prp:eq_gengeo} implies that $S_G$ is an sg-set of $G$. Let $\widetilde{I}(S)$ be a set of paths corresponding to the sg-set $S$ and let $\widetilde{I}(S_G) = p_G(\widetilde{I}(S))$. 

Let $H_1$ and $H_2$ be convex subgraphs of $H$ induced by a convex $2$-partition. We partition $S$ into sets  $S_1$ and $S_2$ such that $S_1 = \{(u,v) \in S:\ v\in H_1\}$ and $S_2 = \{(u,v) \in S:\ v\in H_2\}$. We can assume without loss of generality that $|S_1| \ge |S_2|$. Hence $|S_2| \le \sg(G)/2$. 

We claim that in $\widetilde{I}(S_G)$ the paths between the vertices of $p_G(S_2)$ and the vertices of $S_G$ are not sufficient to cover $G$. Indeed, if this would be the case, then $p_G(S_2)$ would be a strong geodetic core of size at most $\sg(G)/2$, but this is not possible because $\sgc(G) > \sg(G)/2$. By the claim, there exist vertices $u,v\in p_G(S_1)$ such that the $u,v$-path from $\widetilde{I}(S_G)$ covers 
a vertex $w$ which is not covered by the paths with at least one end-vertex in $p_G(S_2)$.

Let $t$ be a vertex of $H_2$. Then the vertex $(w,t)$ is not covered with the paths of  $\widetilde{I}(S)$ that have at least one end-vertex in $S_2$, for otherwise $w$ would be covered in $G$ with a path having one end-vertex in $p_G(S_2)$. Since $H_1$ is a convex subgraph of $H$, there is no shortest path between vertices of $H_1$ containing $t$, so there is no shortest path between vertices of $S_1$ containing $(w,t)$. Hence the vertex $(w,t)$ of $G \cp H$ is not covered by shortest paths between vertices of $S$, a contradiction with the fact that $S$ is an sg-set of $G \cp H$.
\qed

Note that if  $H$ contains a simplicial vertex $u$, then it admits a convex $2$-partition, where $\{u\}$ is one part of it. Hence Theorem~\ref{thm:sgc_sup_sg2} holds  for all such graphs $H$. 

\medskip
In the rest of the section we determine an infinite class of graphs for which the equality is attained in Proposition~\ref{prp:gen_geodesic}. 

\begin{figure}[!ht]
	\begin{center}
		\begin{tikzpicture}[scale=0.6]
		\begin{scope}
		\node[noeud, label=below: $l_1$] (1) at (0,0) {};
		\node[noeud, fill=gray] (3) at (1.5,0) {};
		\node[noeud, fill=gray] (4) at (3,0) {};
		\node[noeud, label=above: $l_2$] (5) at (4.5,1) {};
		\node[noeud, label=below: $l_3$] (6) at (4.5,-1) {};
		
		\draw (1) -- (3) -- (4) --(5);
		\draw (4) -- (6);
		\node at (1,1.5) {\large $T$};
		\end{scope}
		
		\begin{scope}[xshift=10cm]
		\node[noeud] (1) at (0,0) {};
		\node[noeud] (2) at (-1,1) {};
		\node[noeud] (0) at (-1,-1) {};
		\node[noeud, fill=gray] (3) at (1.5,0) {};
		\node[noeud, fill=gray] (4) at (3,0) {};
		\node[noeud] (5) at (4.5,1) {};
		\node[noeud] (6) at (4.5,-1) {};
		\node[noeud] (5') at (3.5,1.7) {};
		\node[noeud] (6') at (3.5,-1.7) {};
		
		\draw (2) -- (3) -- (0) -- (2) -- (1) -- (0);
		\draw (1) -- (3) -- (4) --(5);
		\draw (4) -- (6);
		\draw (4) -- (5') -- (5);
		\draw (4) -- (6') -- (6);
		
		\draw[rotate=0, thick] (-0.6,0) ellipse (1cm and 2cm);
		\node at (-2.3,0) {$K_3$};
		
		\begin{scope}[xshift=3.9cm, yshift=1.5cm]
		\draw[rotate=50, thick] (0,0) ellipse (0.7cm and 1.5cm);
		\end{scope}
		\node at (4.8,2.5) {$K_2$};
		
		\begin{scope}[xshift=3.9cm, yshift=-1.5cm]
		\draw[rotate=-50, thick] (0,0) ellipse (0.7cm and 1.5cm);
		\end{scope}
		\node at (4.8,-2.5) {$K_2$};
		\end{scope}
		\end{tikzpicture}
		\caption{A tree $T$ and the clique tree $K^T_{3,2,2}$.}
		\label{fig:clique_tree}
	\end{center}
\end{figure}
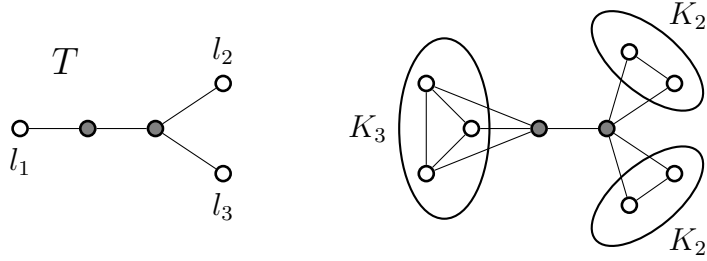

Let $T$ be a tree of order at least $3$ with leaves $l_1, \ldots, l_{\ell(T)}$ and let $n_i$, $i\in [\ell(T)]$, be positive integers. Then the {\em clique tree} $K^{T}_{n_1, \ldots, n_{\ell(T)}}$  of $T$, is the graph obtained from $T$ is the following way. Each leaf $l_i$ of $T$ is replaced with a complete graph of order $n_i$ and each vertex of the clique is connected with an edge to the support vertex of $l_i$. For the case $T=K_2$ we set $K^{K_2}_{n_1,n_2} = K_{n_1+n_2}$. Note that $K^{T}_{1, \ldots, 1} = T$. See Fig.~\ref{fig:clique_tree} for a tree $T$ and its clique tree $K^{T}_{3,2,2}$, where the support vertices of the three leaves are drawn gray. 

%The set of simplicial vertices of the clique tree $K^{T}_{n_1, \ldots, n_{\ell(T)}}$ is 
%$$S = \bigcup_{i \in [\ell(T)]} V(K_{n_i}) \, ,$$
%and is also a strong geodetic set. Hence,
%$$\sg(\CT_{T; n_1, \ldots, n_{\ell(T)}}) = \sum_{i = 1}^{\ell(T)} \max\{1, n_i\} \, .$$ Notice also that $\sgc(\CT_{T; n_1, \ldots, n_{\ell(T)}}) = 1$.

\begin{theorem}
	\label{thm:clique_trees}
	Let $T$ be a tree, let $n_1, \ldots, n_{\ell(T)}$ be positive integers, and let $s = \sg(K^{T}_{n_1, \ldots, n_{\ell(T)}})$. If $G$ is a connected graph with $n(G) \leq s/2$, then 
	$$\sg(K^{T}_{n_1, \ldots, n_{\ell(T)}} \cp G) = s \, .$$
\end{theorem}

\proof
Set $K = K^{T}_{n_1, \ldots, n_{\ell(T)}}$, $n = |V(G)|$, and $V(G) = [n]$. Thus the $K$-layers of the product $K \cp G$ will be denoted by $K^1, \ldots, K^n$. 

Note first that the set of simplicial vertices of $K$ is 
$S = \bigcup_{i \in [\ell(T)]} V(K_{n_i})$, and that $S$ is also a strong geodetic set. Hence, $\sg(K) = \sum_{i = 1}^{\ell(T)} n_i$. Observe also that $\sgc(K) = 1$ and that, due to the tree structure of $K$, the clique tree is geodetic and every simplicial vertex is a core. 

Select and fix $2n$ distinct vertices in $S$: $a_1, b_1, \ldots, a_n, b_n$. Note that this is possible, because $n \leq s/2$. Let $S' = S - \{ a_1, b_1, \ldots, a_n, b_n \}$ and let $$S_{K \cp G} = \{(a_i, i), (b_i, i) : \ i \in [n] \} \cup \{(u,1) : \ u \in S'\} \, .$$ In the following we prove that $S_{K \cp G}$ is a strong geodetic set of $K \cp G$. As $|S_{K \cp G}| = 2n + (s - 2n) = s$ and $\sg(K \cp G) \geq s$ (by Proposition~\ref{prp:gen_geodesic}), this concludes the proof.

Denote a shortest path between vertices $x, y$ in the $H$-layer as $x \rightsquigarrow_H y$ for $H \in \{K, G\}$. Recall that the shortest paths in $K$ are unique because $K$ is geodetic. Now fix the following geodesics, for $i,j \in [n]$:
$$(a_i, i) \rightsquigarrow_K (s, i) \rightsquigarrow_G (u,1), \ u \in S' \, ,$$
$$(a_i, i) \rightsquigarrow_K (b_i, i)\, ,$$
$$(a_i, i) \rightsquigarrow_K (a_j, i) \rightsquigarrow_G (a_j, j), \ j < i \, ,$$
$$(a_i, i) \rightsquigarrow_K (b_j, i) \rightsquigarrow_G (b_j, j), \ j > i \, ,$$
$$(b_i, i) \rightsquigarrow_K (a_j, i) \rightsquigarrow_G (a_j,j), \ j > i \, ,$$
$$(b_i, i) \rightsquigarrow_K (b_j, i) \rightsquigarrow_G (b_j,j), \ j < i \, .$$
First notice that for each pair of vertices in $S_{K \cp G}$ at most one geodesic is selected. Moreover, as each $a_i$ is a strong geodetic core of $K$, the geodesics having $a_i$ as an end-vertex cover $K$. Hence the only possibly uncovered vertices in each $K^i$-layer lie on the $a_i, a_j$-geodesics for all $j > i$, and those on $a_i, b_j$-geodesics for all $j < i$. Here we omit the notation $(x, i)$, as we are only considering the $K^i$-layer. Next we explain that the vertices on these geodesics are in fact also covered. First suppose $j > i$. Consider the subgraph of $K$ induced on the geodesics between vertices $a_i, a_j, b_i$. As $b_i$ is a strong geodetic core of this subgraph, the $a_i, a_j$-geodesic is covered by the geodesics with the end-vertex $b_i$. Next, suppose $j < i $. By a similar reasoning, geodesics containing $a_i$ cover all vertices, except perhaps those on the $a_i, b_j$-geodesic, which is covered by geodesics containing $b_i$. Hence, each $K^i$-layer is completely covered by the selected geodesics.
\qed

Note that the assumption $n(G) \leq s/2$ is necessary. Consider for example the prism over a path $P_n$, $n \geq 2$. Here, $n(K_2) > 2/2$, and $\sg(P_n \cp K_2) = 3 > \sg(P_n)$. 

\begin{corollary}
	\label{cor:trees}
	(i) If $T$ is a tree and $G$ is a connected graph with $2 n(G) \leq \ell(T)$, then $\sg(T \cp G) = \sg(T) = \ell(T)$.
	
	(ii) If $k \geq 2$ and $G$ is a connected graph with $2 n(G) \leq k$, then $\sg(K_k \cp G) = \sg(K_k) = k$.
\end{corollary}

%%%%%%%%%%%%%%%%%%%%%%%%%%%%%%%%%%%%%%%%%%%%%%%%%%%%%%%%%%%%%%%%%%%%%5
%%%%%%%%%%%%%%%%%%%%%%%%%%%%%%%%%%%%%%%%%%%%%%%%%%%%%%%%%%%%%%%%%%%%%5
\section{Counterexamples to Conjecture~\ref{conj:lower-bound}}
\label{sec:counter_conjecture}
%%%%%%%%%%%%%%%%%%%%%%%%%%%%%%%%%%%%%%%%%%%%%%%%%%%%%%%%%%%%%%%%%%%%%5
%%%%%%%%%%%%%%%%%%%%%%%%%%%%%%%%%%%%%%%%%%%%%%%%%%%%%%%%%%%%%%%%%%%%%5

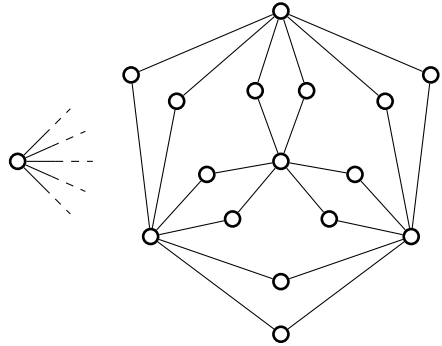
\begin{figure}[!ht]
\begin{center}
\begin{tikzpicture}

\node[noeud] (x) at (0,0){};
\node[noeud] (a) at (90:2){};
\node[noeud] (b) at (-30:2){};
\node[noeud] (c) at (210:2){};

\node[noeud] (xa1) at (70:1){};
\node[noeud] (xa2) at (110:1){};

\node[noeud] (xb1) at (-50:1){};
\node[noeud] (xb2) at (-10:1){};

\node[noeud] (xc1) at (190:1){};
\node[noeud] (xc2) at (230:1){};

\node[noeud] (ab1) at (30:1.6){};
\node[noeud] (ab2) at (30:2.3){};

\node[noeud] (ac1) at (150:1.6){};
\node[noeud] (ac2) at (150:2.3){};

\node[noeud] (bc1) at (270:1.6){};
\node[noeud] (bc2) at (270:2.3){};

\draw (x) -- (xa1);
\draw (x) -- (xa2);
\draw (x) -- (xb1);
\draw (x) -- (xb2);
\draw (x) -- (xc1);
\draw (x) -- (xc2);

\draw (a) -- (xa1);
\draw (a) -- (xa2);
\draw (a) -- (ab1);
\draw (a) -- (ab2);
\draw (a) -- (ac1);
\draw (a) -- (ac2);

\draw (b) -- (ab1);
\draw (b) -- (ab2);
\draw (b) -- (xb1);
\draw (b) -- (xb2);
\draw (b) -- (bc1);
\draw (b) -- (bc2);

\draw (c) -- (xc1);
\draw (c) -- (xc2);
\draw (c) -- (bc1);
\draw (c) -- (bc2);
\draw (c) -- (ac1);
\draw (c) -- (ac2);

\node[noeud] (u) at (-3.5,0){};
\draw (u) -- (-3,0);
\draw[dashed] (-3,0) -- (-2.5,0);
\draw (u) -- (-3.15,0.35);
\draw[dashed] (-3.15,0.35) -- (-2.8,0.7);
\draw (u) -- (-3.05,0.2);
\draw[dashed] (-3.05,0.2) -- (-2.6,0.4);
\draw (u) -- (-3.15,-0.35);
\draw[dashed] (-3.15,-0.35) -- (-2.8,-0.7);
\draw (u) -- (-3.05,-0.2);
\draw[dashed] (-3.05,-0.2) -- (-2.6,-0.4);

\end{tikzpicture}
\end{center}
\caption{A representation of $G_{4,2}$}
\label{fig:g42}
\end{figure}

In this section we demonstrate that, rather surprisingly, Conjecture~\ref{conj:lower-bound} does not hold in general. For this sake consider the graphs $G_{k,n}$, $k \geq 4$, $n \geq 2$, constructed as follows. Start with the complete graph $K_k$, replace every edge of it with $n$ disjoint paths of length $2$, and finally add an additional universal vertex, that is, a vertex adjacent to all other vertices. In Fig.~\ref{fig:g42} the graph $G_{4,2}$ is shown. 

\begin{theorem}
\label{thm:counterexample}
If $k \geq 4$ and $n \geq 2$, then 
$$\sg(G_{k,n}) = \binom{k}{2} (n-1) + k \quad \text{and} \quad \sg(G_{k,n} \cp K_n) \leq k n + 1\,.$$
\end{theorem}

\proof 
Let $k \geq 4$ and $n \geq 2$ and denote the vertices of $G_{k,n}$ as follows. Let $x_1,\ldots, x_k$ be its vertices corresponding to the initial complete graph $K_k$, let $x_{ij}^{(1)}, \ldots, x_{ij}^{(n)}$ be the common degree-$2$ neighbors of $x_i$ and $x_j$, where $ i,j \in [k]$, $i \neq j$, and let $u$ be the universal vertex of $G_{k,n}$. 

An arbitrary degree-2 vertex $x_{ij}^{(l)}$ of $G_{k,n}$ can be covered by either being in a strong geodetic set, or by the unique $x_i,x_j$-geodesic. Thus, for all $\{i, j\}\in \binom{[k]}{2}$, at least $n-1$ vertices among $x_{ij}^{(1)}, \ldots, x_{ij}^{(n)}$ must lie in any strong geodetic set. 

Notice that $S = \{x_1,\ldots, x_k\} \cup \{ x_{ij}^{(l)}:\ i, j \in [k], i \neq j, l \in [n-1] \}$ is a strong geodetic set of size $k + \binom{k}{2} (n-1)$. Suppose there exists a smaller strong geodetic set $T$. We may without loss of generality assume that $x_1 \notin T$. Then $T$ contains vertices $x_{1j}^{(l)}$ for $j \in \{2,\ldots, k\}$ and $l \in [n]$. If also some other $x_i \notin T$, then $|T| \geq \binom{k}{2} (n-1) + (k-1) + (k-2) > |S|$ (as $k \geq 4$). But if all $x_2, \ldots, x_k \in T$, then $|T| \geq \binom{k}{2} (n-1) + (k-1) + (k-1) > |S|$. Hence, $\sg(G_{k,n}) = |S|$. 

Consider now the Cartesian product graph $G_{k,n} \cp K_n$. Set $V(K_n) = \{y_1, \ldots, y_n\}$ and 
$$S = \{ (u, y_1) \} \cup \{ (x_i, y_j):\ i \in [k], j \in [n] \}\,.$$ 
Fix geodesics in the following way
\begin{align*}
	(x_i, y_j) \sim (x_{ii'}^{(j')}, y_j) \sim (x_{i'}, y_j) \sim (x_{i'}, y_{j'}); & \quad  i, i' \in [k], i < i', j, j' \in [n], j \neq j',\\
	(x_i, y_j) \sim (x_{ii'}^{(j)}, y_j) \sim (x_{i'}, y_j); & \quad  i, i' \in [k], i < i', j \in [n],\\
	(u, y_1) \sim (u, y_j) \sim (x_1, y_j); & \quad j \in \{ 2, \ldots, n \},
\end{align*}
to see that $S$ is a strong geodetic set of $G_{k,n} \cp K_n$. We conclude that $\sg(G_{k,n} \cp K_n) \leq |S| = 1 + k n$. 
\qed

It follows from Theorem~\ref{thm:counterexample} that for every $k \geq 4$ and every  $n \geq 2$, 
$$\sg(G_{k,n}) - \sg(G_{k,n} \cp K_n) \ge \frac{k(n-1)(k-3) - 2}{2}\,.$$ 
Hence by increasing $k$ or $n$, the difference $\sg(G_{k,n}) - \sg(G_{k,n} \cp K_n)$ becomes arbitrary large. This thus disproves Conjecture~\ref{conj:lower-bound}. 

\medskip
Note that the counterexamples from Theorem~\ref{thm:counterexample} are of diameter $2$. For counterexamples of higher diameter consider the following construction. Let $H$ be the graph obtained from a complete graph $K_k$ in which every edge is replaced with $n$ pairwise disjoint paths of length $2p$, $p\ge 2$. Let us call the middle vertex of such a path in $H$ a {\em middle vertex}, and let the vertices of $H$ that correspond to the vertices of $K_k$ be called {\em original vertices}. Next, let $H'$ be the graph obtained from the star $K_{1,k}$ by replacing each of its edges with a path of length $p-1$. Finally, connect each leaf of $H'$ with a unique original vertex of $H$, and connect the vertex of $H'$ of degree $k$ with all the middle vertices of $H$. Then the constructed graph has diameter $2p$ and by similar yet more technical arguments as in the proof of Theorem~\ref{thm:counterexample} one can show that the Cartesian product of it with $K_n$ is also a counterexample.   

%%%%%%%%%%%%%%%%%%%%%%%%%%%%%%%%%%%%%%%%%%%%%%%%%%%%%%%%%%%%%%%%%%%%%5
%%%%%%%%%%%%%%%%%%%%%%%%%%%%%%%%%%%%%%%%%%%%%%%%%%%%%%%%%%%%%%%%%%%%%5
\section{Concluding remarks}
\label{sec:concluding}
%%%%%%%%%%%%%%%%%%%%%%%%%%%%%%%%%%%%%%%%%%%%%%%%%%%%%%%%%%%%%%%%%%%%%5
%%%%%%%%%%%%%%%%%%%%%%%%%%%%%%%%%%%%%%%%%%%%%%%%%%%%%%%%%%%%%%%%%%%%%5

In this paper we have studied the strong geodetic problem on graphs. The newly introduced strong geodetic cores turned out to be quite useful for this sake.   Hence we believe that to introduce and investigate {\em geodetic cores} for the usual geodetic problem would be a reasonable research plan.  

We have also introduced generalized geodetic graphs as the graphs $G$ for which $\g(G) = \sg(G)$ holds. In Section~\ref{sec:valid_conjecture} we showed that these graphs behave nicely on the Cartesian product w.r.t.\ the strong geodetic problem. Hence it seems an interesting problem to characterize the class of generalized geodetic graphs. 

It would also be interesting to characterize the graphs which have unique sg-sets, because for them it seems easier to determine strong geodetic cores than in the general case.  

Finally, the computational complexity of the strong geodetic core number would be interesting to investigate.  

\section*{Acknowledgements}

Research supported from the Slovenian Research Agency (research core funding No. P1-0297). 

%%%%%%%%%%%%%%%%%%%%%%%%%%%%{%%%%%%%%%%%%%%%%%%%%%%%%%%%%%%%%%%%%%%

\end{document}